\documentclass[11pt]{article}

\usepackage{latexsym}
\usepackage[dvips]{graphicx}
\usepackage{float}
\usepackage{amsmath}
\usepackage{epic}
\usepackage{color}
\usepackage{setspace}
\usepackage{lscape}
\usepackage{arydshln}

\begin{document}
\bibliographystyle{plain}
\floatplacement{table}{H}
\newtheorem{definition}{Definition}[section]
\newtheorem{lemma}{Lemma}[section]
\newtheorem{theorem}{Theorem}[section]
\newtheorem{corollary}{Corollary}[section]
\newtheorem{proposition}{Proposition}[section]
\newcommand{\sni}{\sum_{i=1}^{n}}
\newcommand{\snj}{\sum_{j=1}^{n}}
\newcommand{\smj}{\sum_{j=1}^{m}}
\newcommand{\sumjm}{\sum_{j=1}^{m}}
\newcommand{\bdis}{\begin{displaymath}}
\newcommand{\edis}{\end{displaymath}}
\newcommand{\beq}{\begin{equation}}
\newcommand{\eeq}{\end{equation}}
\newcommand{\beqn}{\begin{eqnarray}}
\newcommand{\eeqn}{\end{eqnarray}}
\newcommand{\qed}{{\large $\sqcap$ \hskip -0.37cm $\sqcup$}}
\newcommand{\defeq}{\stackrel{\triangle}{=}}
\newcommand{\sep}{\;\;\;\;\;\; ; \;\;\;\;\;\;}
\newcommand{\real}{\mbox{$ I \hskip -4.0pt R $}}
\newcommand{\complex}{\mbox{$ I \hskip -6.8pt C $}}
\newcommand{\integ}{\mbox{$ Z $}}
\newcommand{\realn}{\real ^{n}}
\newcommand{\sqrn}{\sqrt{n}}
\newcommand{\sqrtwo}{\sqrt{2}}
\newcommand{\prf}{{\bf Proof. }}
\newcommand{\onehlf}{\frac{1}{2}}
\newcommand{\thrhlf}{\frac{3}{2}}
\newcommand{\fivhlf}{\frac{5}{2}}
\newcommand{\onethd}{\frac{1}{3}}
\newcommand{\lb}{\left ( }
\newcommand{\lcb}{\left \{ }
\newcommand{\lsb}{\left [ }
\newcommand{\labs}{\left | }
\newcommand{\rb}{\right ) }
\newcommand{\rcb}{\right \} }
\newcommand{\rsb}{\right ] }
\newcommand{\rabs}{\right | }
\newcommand{\lnm}{\left \| }
\newcommand{\rnm}{\right \| }
\newcommand{\lambdab}{\bar{\lambda}}
%
%

\newcommand{\sumjrosig}{\sum_{\stackrel{j=1}{j\neq\rho,\sigma}}^{n}}
\newcommand{\summuro}{\sum_{\stackrel{j=1}{j\neq\mu,\rho}}^{n}}
\newcommand{\sumjnoti}{\sum_{\stackrel{j=1}{j\neq i}}^{n}}
\newcommand{\sumlnoti}{\sum_{\stackrel{\ell=1}{\ell \neq i}}^{n}}
\newcommand{\sumknoti}{\sum_{\stackrel{k=1}{k\neq i}}^{n}}
\newcommand{\sumknotij}{\sum_{\stackrel{k=1}{k\neq i,j}}^{n}}
\newcommand{\sumk}{\sum_{k=1}^{n}}
\newcommand{\snl}{\sum_{\ell=1}^{n}}

\newcommand{\xj}{x_{j}}
\newcommand{\xjpo}{x_{j+1}}
\newcommand{\xjmo}{x_{j-1}}

\newcommand{\azero}{a_{0}}
\newcommand{\aone}{a_{1}}
\newcommand{\atwo}{a_{2}}
\newcommand{\ath}{a_{3}}
\newcommand{\afr}{a_{4}}
\newcommand{\afv}{a_{5}}
\newcommand{\asx}{a_{6}}
\newcommand{\anmo}{a_{n-1}}
\newcommand{\akmo}{a_{k-1}}
\newcommand{\akpo}{a_{k+1}}
\newcommand{\anmt}{a_{n-2}}
\newcommand{\akmt}{a_{k-2}}
\newcommand{\an}{a_{n}}
\newcommand{\ak}{a_{k}}

\newcommand{\bzero}{b_{0}}
\newcommand{\bone}{b_{1}}
\newcommand{\btwo}{b_{2}}
\newcommand{\bth}{b_{3}}
\newcommand{\bfr}{b_{4}}
\newcommand{\bfv}{b_{5}}
\newcommand{\bsx}{b_{6}}
\newcommand{\bnmo}{b_{n-1}}
\newcommand{\bkmo}{b_{k-1}}
\newcommand{\bkpo}{b_{k+1}}
\newcommand{\bnmt}{b_{n-2}}
\newcommand{\bkmt}{b_{k-2}}
\newcommand{\bn}{b_{n}}
\newcommand{\bk}{b_{k}}

\newcommand{\alfazero}{\alpha_{0}}
\newcommand{\alfaone}{\alpha_{1}}
\newcommand{\alfatwo}{\alpha_{2}}
\newcommand{\alfath}{\alpha_{3}}
\newcommand{\alfafr}{\alpha_{4}}
\newcommand{\alfafv}{\alpha_{5}}
\newcommand{\alfasx}{\alpha_{6}}
\newcommand{\alfanmo}{\alpha_{n-1}}
\newcommand{\alfanmt}{\alpha_{n-2}}
\newcommand{\alfan}{\alpha_{n}}
\newcommand{\alfak}{\alpha_{k}}
\newcommand{\bkl}{b_{k\ell}}

\newcommand{\xb}{\bar{x}}        
\newcommand{\wb}{\bar{w}}        
\newcommand{\yb}{\bar{y}}        
\newcommand{\ub}{\bar{u}}        
\newcommand{\phix}{\phi(x)}      
\newcommand{\phixb}{\phi(\xb)}      
\newcommand{\phione}{\phi_{1}}      
\newcommand{\phitwo}{\phi_{2}}      
\newcommand{\phionex}{\phi_{1}(x)}      
\newcommand{\phionexb}{\phi_{1}(\xb)}      
\newcommand{\phitwox}{\phi_{2}(x)}      
\newcommand{\phitwoxb}{\phi_{2}(\xb)}      

\newcommand{\rone}{r}
\newcommand{\rtwo}{R}
\newcommand{\obar}{\widebar{O}}
\newcommand{\Lbar}{\widebar{L}}
\newcommand{\deltaone}{\delta_{1}}
\newcommand{\deltatwo}{\delta_{2}}
\newcommand{\rhoone}{\rho_{1}}
\newcommand{\rhotwo}{\rho_{2}}

\newcommand{\thzero}{\eta_{0}}
\newcommand{\thone}{\eta_{1}}
\newcommand{\thnmo}{\eta_{n-1}}
\newcommand{\thnmt}{\eta_{n-2}}
\newcommand{\thk}{\eta_{k}}
\newcommand{\thj}{\eta_{j}}
\newcommand{\thkmo}{\eta_{k-1}}
\newcommand{\thkmt}{\eta_{k-2}}
\newcommand{\thkpo}{\eta_{k+1}}

\begin{center}
\large
{\bf IMPLEMENTATION OF PELLET'S THEOREM} \\
\vskip 0.5cm
\normalsize
A. Melman \\
Department of Applied Mathematics \\
School of Engineering, Santa Clara University  \\
Santa Clara, CA 95053  \\
e-mail : amelman@scu.edu \\
\vskip 0.5cm
\end{center}
\vskip 0.2cm
\begin{abstract}
Pellet's theorem determines when the zeros of a polynomial can be separated into two regions,
based on the presence or absence of positive roots of an auxiliary polynomial, but does not 
provide a method to verify its conditions or to compute the roots of the auxiliary polynomial 
when they exist. We derive an explicit condition
for these roots to exist and, when they do, propose efficient ways to compute them.
A similar auxiliary polynomial appears for the generalized Pellet theorem for matrix polynomials
and it can be treated in the same way.

\vskip 0.15cm
{\bf Key words :} Pellet, zero, root, polynomial, matrix polynomial
\vskip 0.15cm
{\bf AMS(MOS) subject classification :} 12D10, 30C15
\end{abstract}


\section{Introduction}               


Pellet's theorem, a classical result from 1881, derives conditions under which the zeros of a polynomial can be divided
into two groups, according to their magnitudes. It is a direct consequence of Rouch\'{e}'s
theorem and is stated as follows.
\begin{theorem}
(\cite{Pellet}, \cite[Th.(28,1), p.128]{Marden})
\label{Pellet}
Given the polynomial $p(z)=z^{n} + a_{n-1} z^{n-1} + \dots + \aone z + \azero$ with complex coefficients, 
$n \geq 3$, $1 \leq k \leq n-1$, and $\azero \ak \neq 0$. Let the polynomial 
\bdis
\phi(x)=x^{n} + |a_{n-1}| x^{n-1} + \dots + |a_{k+1}| x^{k+1} - |a_{k}| x^{k} + |a_{k-1}|x^{k-1} + \dots  + |\azero| 
\edis
have two distinct positive roots $r$ and $R$, $r < R$.
Then $p$ has exactly $k$ zeros in or on the circle $|z| = r$ and no zeros in the annular ring $r < |z| < R$.
\end{theorem}
We note that, by Descartes' rule, $\phi$ has either two or no positive roots.
Although the function $\phi$ depends on $k$, this parameter
remains fixed throughout the paper, and it will therefore be omitted from the notation to prevent unnecessary clutter.

Recently, (\cite{BiniNoferiniSharify}, \cite{MelmanMatPol}) a generalized
Pellet theorem was derived for matrix polynomials, which 
have received a lot of attention recently because of their application in several 
engineering fields (\cite{TisseurMeerbergen}). Matrix polynomials are polynomials  
whose coefficients are matrices instead of scalars. They occur in polynomial eigenvalue 
problems, which consist of finding a nonzero eigenvector $v$, corresponding to an eigenvalue $z$
satisfying $P(z)v=0$, where
\bdis
P(z) = A_{n} z^{n} + A_{n-1} z^{n-1} + \dots + A_{0},
\edis
with $A_{j} \in \complex^{m \times m}$ for $j=0,\dots n$, and with $\det{(P(z))}$ not
identically zero.
If $A_{n}$ is singular then $P$ has infinite eigenvalues and if $A_{0}$ is singular then zero 
is an eigenvalue.  There are $nm$ eigenvalues, including possibly infinite ones. The finite 
eigenvalues are the solutions of $\det{(P(z))}=0$.
The following theorem generalizes Pellet's theorem to matrix polynomials.
%
%
\begin{theorem}
\label{genPellet}             
{\bf (Generalized Pellet theorem.)} (\cite{BiniNoferiniSharify}, \cite{MelmanMatPol})
Let 
\bdis  
\label{Pelleteq}
P(z) = A_{n}z^{n} + A_{n-1}z^{n-1} + \dots + A_{1}z + A_{0}  \nonumber 
\edis
be a matrix polynomial with $n \geq 2$, $A_{j} \in \complex^{m \times m}$ for $j=0,\dots,n$, and $A_{0} \neq 0$.      
Let $A_{k}$ be invertible for some $k$ with $1 \leq k \leq n-1$, and let the polynomial 
\bdis
g(x)=||A_{n}||x^{n} + ||A_{n-1}||x^{n-1} + \dots + ||A_{k+1}||x^{k+1} - ||A^{-1}_{k}||^{-1}x^{k} 
+ ||A_{k-1}||x^{k-1} + \dots + ||A_{1}||x + ||A_{0}|| \; ,   
\edis
where the norm can be any vector-induced norm, have two distinct positive roots $r$ and $R$, $r < R$.
Then $P$ has exactly $km$ zeros in or on the disk $|z|=r$ and no zeros in the annular ring 
$r < |z| < R$.
\end{theorem}
Its generalization widens the usefulness of Pellet's theorem considerably, as the computation of nonlinear
eigenvalues is much more costly than the computation of polynomial zeros, making easily computed bounds more valuable.

However, both the classical and generalized Pellet theorems do not explain how their assumptions can be verified, 
nor do they provide a method to compute the roots $r$ and $R$ when the assumptions are satisfied. In addition,
as will be explained below, we require certain properties of any iterates generated in such computations. 
Our sole purpose is to address these shortcomings, thereby improving the applicability of these theorems.
We stress that we focus on the implementation of Pellet's theorem, \emph{if and when}
it is used. When and where it is worthwile to use the theorem is a consideration that lies well 
outside the scope of this work.

Because the real polynomials $\phi$ and $g$ are of the exact same form, the techniques we will develop
apply to both Theorem~\ref{Pellet} and Theorem~\ref{genPellet}. To keep matters simple, we will
henceforth refer only to Theorem~\ref{Pellet}, with the understanding that all results carry
over in a straightforward way to Theorem~\ref{genPellet}.

\begin{figure}[H]
\begin{center}
\hskip -0.5cm
\includegraphics[width=0.30\linewidth]{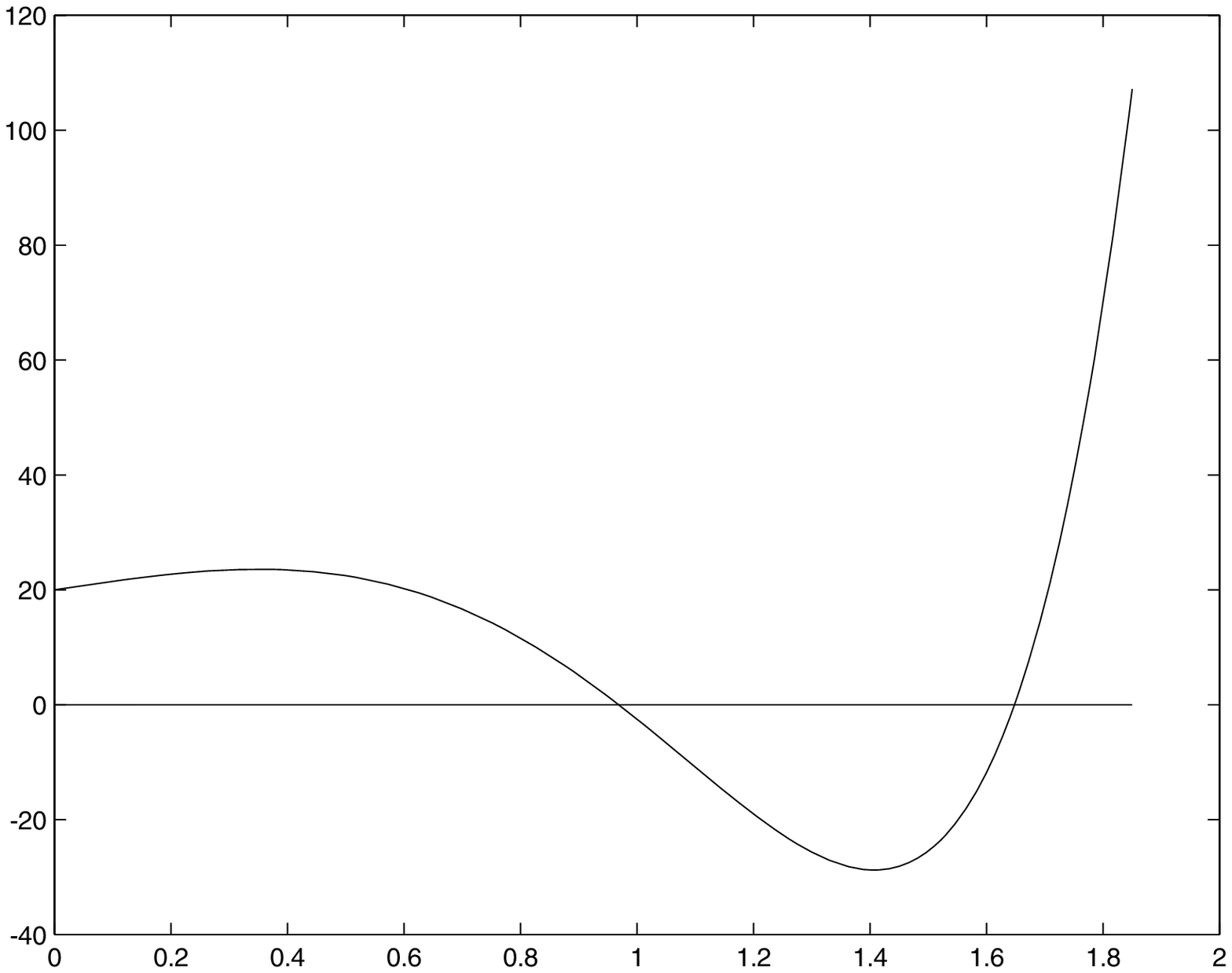}
\;\;   
\includegraphics[width=0.30\linewidth]{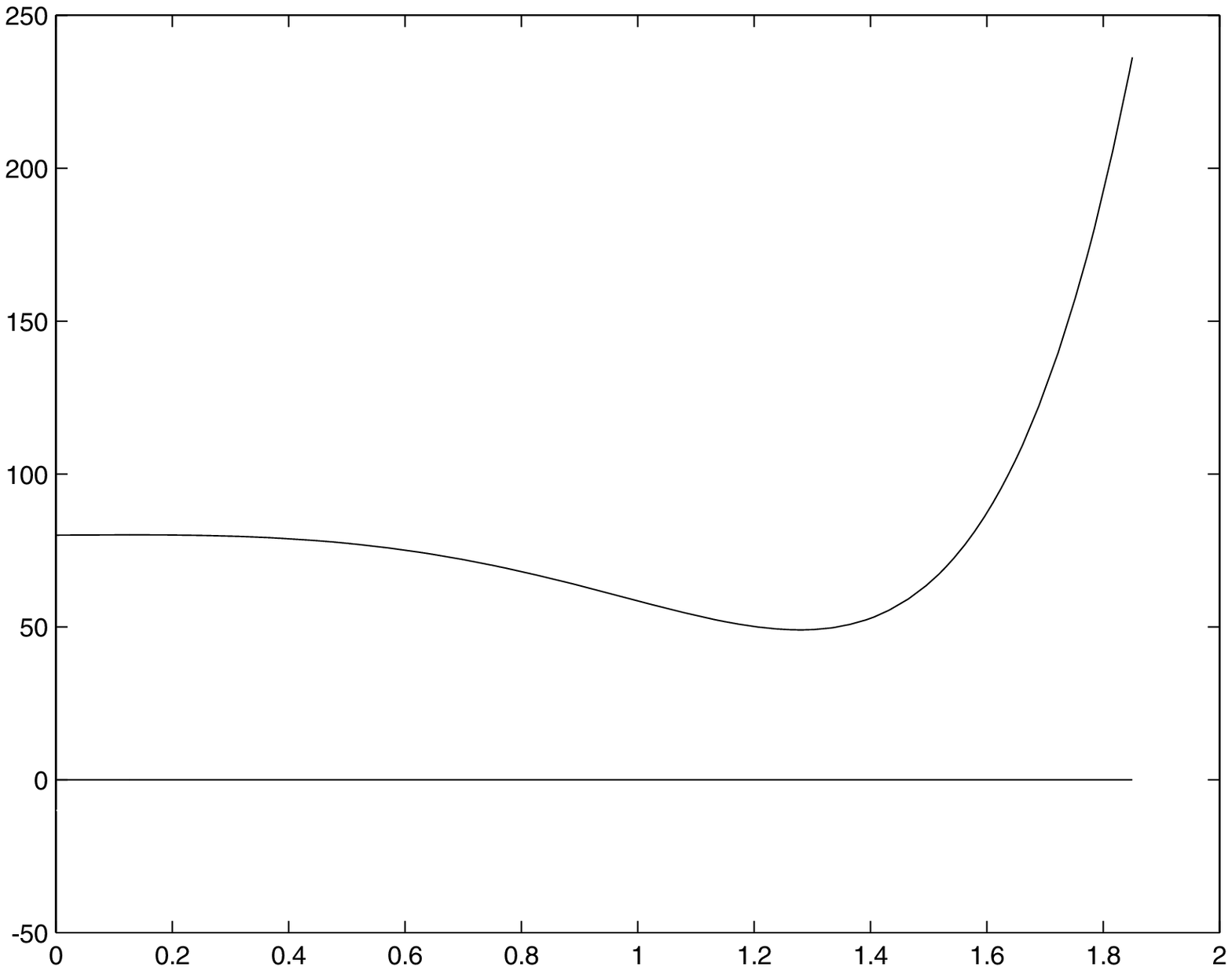}
\;\;
\includegraphics[width=0.30\linewidth]{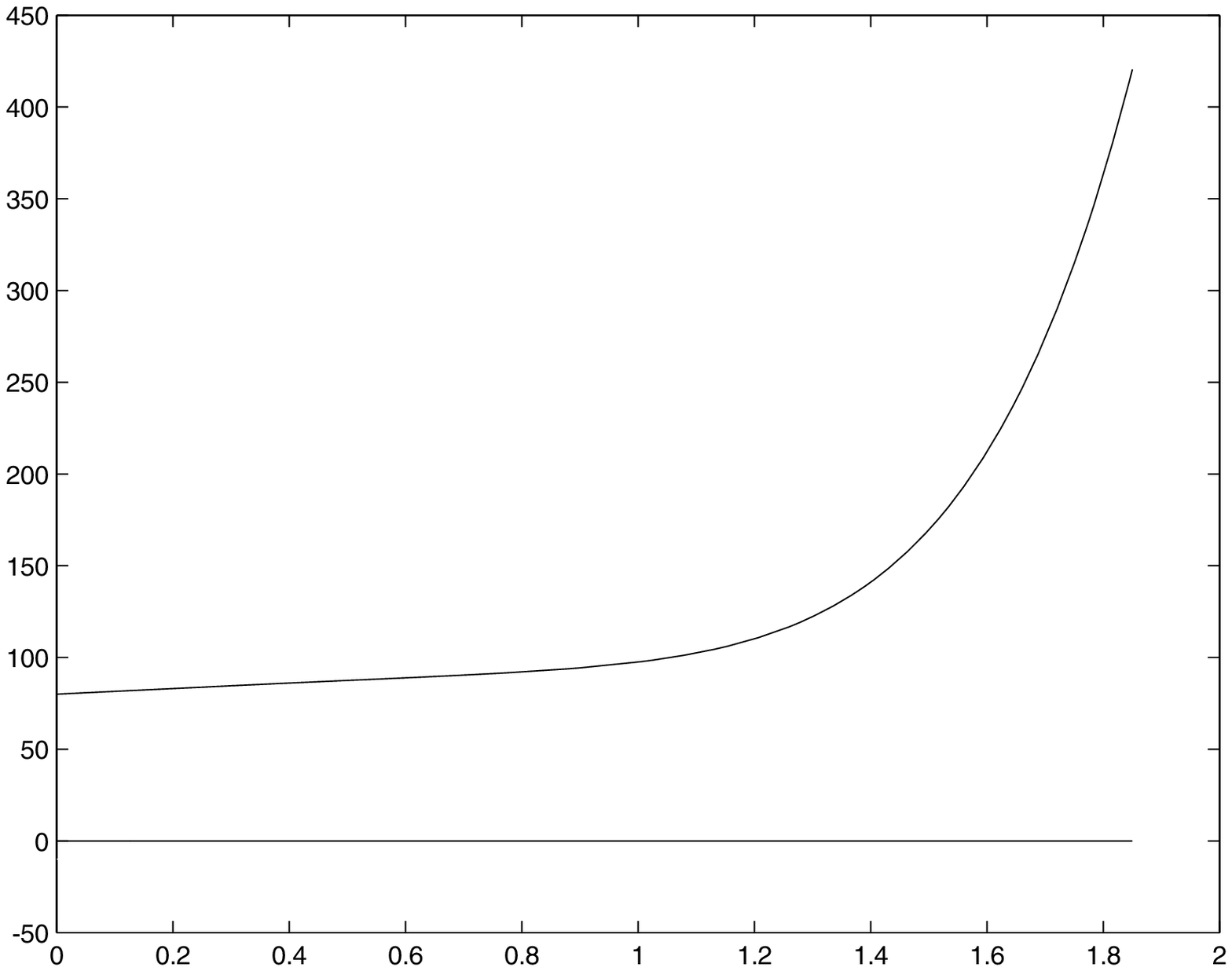}
\end{center}
\caption{Typical scenarios for $\phi$.}
\label{phiscenarios}    
\end{figure}

Figure~\ref{phiscenarios} shows a few typical scenarios for the function $\phi$ in Theorem~\ref{Pellet}:
on the left $\phi$ increases at first, then decreases , intersects the x-axis and then increases while crossing
the x-axis again; in the middle, $\phi$ behaves as on the left, but does not decrease enough to have positive roots; 
on the right $\phi$ increases monotonically.
The presence or lack of positive roots needs to be detected first, and, 
if there are such roots, then they need to be computed with a method that generates iterates   
that are themselves proper bounds, which means that
$r$ and $R$ need to be approached from inside the interval $[r,R]$. 
In this way, the numerical process can be stopped at any time with an upper bound on $r$ and a lower bound on $R$,
allowing for an inexact solution while still providing correct bounds on the two groups of $k$ and $n-k$ zeros~of~$p$.
We consider this a key property to be satisfied and, although there exist many methods to compute the real roots
of a polynomial, none of them accomplishes this.  

We derive an easily computable criterion for the absence or presence of positive roots of $\phi$ in Theorem~\ref{Pellet}, 
which has the important additional advantage of providing an adequate starting point in $[r,R]$ for the computation 
of the positive roots of $\phi$.
This allows us to reformulate Pellet's theorem in a more useful way and to propose a framework for generating 
efficient methods to compute the roots, when they exist, from inside $[r,R]$ for a particular given value of $k$. 
The complexity of detecting whether $\phi$ has positive roots or not is, as we will see later, $\mathcal{O}(n)$. 
This means that, if the theorem needs to be applied for \emph{every} value of $k$,
with $1 \leq k \leq n-1$, then this complexity becomes $\mathcal{O}(n^{2})$, which makes it more efficient to 
first use a result by \cite{BiniNumAlg}. In \cite{BiniNumAlg},
it was shown that if $\phi$ has two real roots for a particular value of $k$, then $k$ must be the abscissa of a vertex
of the Newton polygon associated with $p$. Moreover, computing these vertices only costs $\mathcal{O}(n\log{n})$ operations
(\cite{Graham}).
In such a case, we would first compute the abscissae $k_{1},...k_{m}$ of the vertices of the associated Newton 
polygon and then compute the roots $r$ and $R$ \emph{only} for the values of $k \in \{k_{1},...,k_{m} \}$. 
Once the abscissae have been computed, the complexity of detecting the positive roots is then  
$\mathcal{O}(mn)$. Typically, $m \ll n$ since only polynomials with very special coefficients have more 
than just a few values of $k$ for which Pellet's theorem can be applied. 

From now on, we compute the roots of $\phi$ for a given value of $k$ which remains fixed throughout 
this work.
Detection of the roots and the determination of a starting point is the subject of 
Section~\ref{detection}, while a strategy to compute the roots is developed in Section~\ref{computation}.
We believe that the ideas behind our techniques are general enough to be useful in other situations as well.
%
%


\section{Detection of the roots}
\label{detection}

Our strategy to determine if $\phi$ has positive roots is to first transform it to 
a strictly convex function $\chi$ with the same roots or absence of roots as $\phi$. 
Denoting $\chi$'s unique minimizer by $x^{*}$, we conclude that
$\phi$ has two positive roots if and only if $\phi(x^{*})<0$.
We note that this minimizer does not, in general, correspond to the minimizer of $\phi$ (if it exists).
This leads to the following theorem.

\begin{theorem}
\label{detectiontheorem}
Given the polynomial 
\bdis
\phi(x)=x^{n} + \thnmo x^{n-1} + \dots + \thkpo x^{k+1} - \thk x^{k} + \thkmo x^{k-1} + \dots  + \thzero  \; , 
\edis
with $\thj \geq 0$ ($j=1,...,n$), $\thzero \thk \neq 0$, $n \geq 3$, and $1 \leq k \leq n-1$. 
Then $\phi$ has two positive roots if and only if $\phi(x^{*})<0$, where $x^{*}$ is the unique positive root of the 
polynomial
\bdis
\chi(x) = (n-k)x^{n} + (n-k-1) \thnmo x^{n-1} + \dots + \thkpo x^{k+1} - \thkmo x^{k-1} -2\thkmt x^{k-2}
- \dots - k\thzero  \; .
\edis
\end{theorem}
\prf
We define $\psi(x)$ for $x \in (0,+\infty)$ as 
\bdis
\psi(x) = x^{-k}\phi(x) 
= x^{n-k} + \thnmo x^{n-k-1} + \dots + \thkpo x - \thk + \thkmo x^{-1} + \dots + \thzero x^{-k}  \; ,
\edis
and observe that it is a strictly convex function with the same positive roots or lack thereof as $\phi$.
This function will therefore have a unique positive minimizer, which we denote by $x^{*}$, and it will have two
positive roots if and only if $\psi(x^{*}) < 0$, which is equivalent to $\phi(x^{*}) < 0$.
Because $x^{*}$ is the unique minimizer of a strictly convex differentiable function, it can be obtained as the unique
positive solution of $\psi'(x)=0$. Since
\bdis
\psi'(x) = (n-k) x^{n-k-1} + (n-k-1) \thnmo x^{n-k-2} + \dots + \thkpo  - \thkmo x^{-2} + \dots - k \thzero x^{-k-1} \; ,
\edis
and since $\psi'(x)=0 \Longleftrightarrow x^{k+1}\psi'(x)=0$, the statement of the theorem follows. \qed

Defining
\bdis
\sigma(x)=x^{n} + \thnmo x^{n-1} + \dots + \thkpo x^{k+1} + \thkmo x^{k-1} + \dots  + \thzero  \; , 
\edis
we remark that Theorem~\ref{detectiontheorem} can also be obtained from the observation that the 
existence of $u > 0$ for which $\phi(u) < 0$ is equivalent to the existence of $u > 0$ for which
$|\ak|>\sigma(u)/u^{k}$. 
The function in the right-hand side of this inequality is strictly convex and has a unique minimizer $u^{*}$. Then
$\sigma(u^{*})/(u^{*})^{k}$ represents the lowest strict threshhold value for $|\ak|$ to guarantee the existence of
some $u > 0$ for which $\phi(u)<0$.

Since $x^{*}$ is independent of the value of $\thk$, it is possible to establish
a criterion for $\thk$ that guarantees positive roots for $\phi$, and that depends only on the other coefficients.
Incorporating this in Pellet's theorem results in the following more explicit version of that same theorem. 

\begin{theorem}
\label{newPellet}
Given the polynomial $p(z)=z^{n} + a_{n-1} z^{n-1} + \dots + \aone z + \azero$ with complex coefficients, 
$\azero \ak \neq 0$, and $n \geq 3$. Let $1 \leq k \leq n-1$, 
\begin{eqnarray*}
& & \phi(x)=x^{n} + |a_{n-1}| x^{n-1} + \dots + |a_{k+1}| x^{k+1} - |a_{k}| x^{k} 
+ |a_{k-1}| x^{k-1} + \dots  + |\azero| \; ,  \\
& & \sigma(x)=x^{n} + |a_{n-1}| x^{n-1} + \dots + |a_{k+1}| x^{k+1} + |a_{k-1}| x^{k-1} + \dots  + |\azero| \; , 
\end{eqnarray*}
and let $x^{*}$ be the unique positive root
of the polynomial
\bdis
\chi(x) = (n-k)x^{n} + (n-k-1) |\anmo| x^{n-1} + \dots + |\akpo| x^{k+1} - |\akmo| x^{k-1} -2 |a_{k-2}| x^{k-2}
- \dots - k |\azero|   \; .
\edis
If
\beq
\label{akcriterion}
|\ak| > \dfrac{\sigma(x^{*})}{(x^{*})^{k}} \; ,         
\eeq
then the polynomial $\phi$ has two distinct positive roots $r$ and $R$ with $r < R$, 
and the polynomial $p$ has exactly $k$ zeros in or on the circle 
$|z| = r$ and no zeros in the annular ring $r < |z| < R$.
\end{theorem}

The polynomial $\chi$ appearing in the two previous theorems, is strictly convex to the right of its positive 
root, so that, to compute it, Newton's method or a suitably accelerated version of it can be used with guaranteed 
monotonic convergence from the right of the root. As a starting point, any easily computable upper bound 
on the roots of $\phi$ can be used. Such a bound can be found with the following theorem.

\begin{theorem}
\label{upperboundomega}
Let
\bdis
\omega(x) = \eta_{n} x^{n} + \eta_{n-1} x^{n-1} + \dots + \eta_{m} x^{m} - \eta_{\ell} x^{\ell} - \eta_{\ell-1} x^{\ell-1}
- \dots - \eta_{0} \; ,
\edis
with $\eta_{j} \geq 0$ ($j=1,...,n$), $m > \ell$, $\eta_{0} \eta_{\ell} \eta_{m} \eta_{n} \neq 0$, and $n \geq 3$. 
Then the unique positive root $x^{*}$ of $\omega$ satisfies the following.
$ $ \newline {\bf (1)} If $\omega(1) < 0$, then 
\bdis
1 < x^{*} \leq  
\lb \dfrac{\eta_{\ell} + \eta_{\ell-1} + \dots + \eta_{0}}{\eta_{n} + \eta_{n-1} 
+ \dots + \eta_{m}} \rb ^{\frac{1}{m-\ell}} \; \cdot  
\edis
$ $ \newline {\bf (2)} If $\omega(1) > 0$, then 
\bdis
0 < x^{*} \leq  
\lb \dfrac{\eta_{\ell} + \eta_{\ell-1} + \dots + \eta_{0}}{\eta_{n} + \eta_{n-1} 
+ \dots + \eta_{m}} \rb ^{\frac{1}{n}} < 1 \; .          
\edis
\end{theorem}
\prf
By Descartes' rule of signs, the polynomial $\omega$ has a single positive root $x^{*}$ and is therefore 
negative for $0 \leq x < x^{*}$.
If $\omega(1)=0$, then $x^{*}=1$. If $\omega(1) < 0$, then $x^{*} > 1$, and for $x > 1$ we have that 
\beq
\omega(x) \geq  \lb \eta_{n} + \eta_{n-1} + \dots + \eta_{m} \rb  x^{m} 
- \lb \eta_{\ell} + \eta_{\ell-1} + \dots + \eta_{0} \rb x^{\ell} \nonumber \; ,
\eeq
with the equality holding for $x=1$. 
The unique positive root of the right-hand side then provides an upper bound on the positive root of $\omega$.

If $\omega(1) > 0$, then $x^{*} < 1$, and for $x < 1$ we have that 
\beq
\label{vec1}
\omega(x) \geq  \lb \eta_{n} + \eta_{n-1} + \dots + \eta_{m} \rb  x^{n} 
- \lb \eta_{\ell} + \eta_{\ell-1} + \dots + \eta_{0} \rb  \nonumber \; ,
\eeq
with equality for $x=1$.
The unique positive root of the right-hand side once again provides an upper bound on the positive root of $\omega$.
That it is less than $1$ follows from $\omega(1) > 0$. This concludes the proof. \qed

Applying this theorem to the aforementioned polynomial $\chi$ directly yields the following corollary.

\begin{corollary}
\label{chibounds}
Let
\bdis
\chi(x) = (n-k)x^{n} + (n-k-1) \thnmo x^{n-1} + \dots + \thkpo x^{k+1} - \thkmo x^{k-1} -2 \thkmt x^{k-2}
- \dots - k\thzero \; ,
\edis
with $\thj \geq 0$ ($j=1,...,n$), $\thzero \neq 0$, $n \geq 3$, and $1 \leq k \leq n-1$. 
Then the unique positive root $x^{*}$ of $\chi$ satisfies the following.
$ $ \newline {\bf (1)} If $\chi(1) < 0$, then 
\beq
\label{chibound1}
1 < x^{*} \leq  
\lb \dfrac{\thkmo +2\thkmt + \dots + k\thzero}{(n-k) + (n-k-1) \thnmo + \dots 
+ \thkpo} \rb ^{1/2} \; \cdot  
\eeq
$ $ \newline {\bf (2)} If $\chi(1) > 0$, then 
\beq
\label{chibound2}
0 < x^{*} \leq  
\lb \dfrac{\thkmo +2\thkmt + \dots + k\thzero}{(n-k) + (n-k-1) \thnmo + \dots 
+ \thkpo} \rb ^{1/n} < 1 \; .          
\eeq
\end{corollary}
The bounds can be adjusted like in Theorem~\ref{upperboundomega} if 
$\eta_{k+1}$, $\eta_{k-1}$, or more coefficients vanish.
There are several possible strategies for finding a proper starting point to begin the computation of the roots 
of $\phi$ if they exist: one can periodically compute $sgn(\phi)$ at an iterate,
and use that iterate as a starting point if $sgn(\phi) = -1$,
or one can simply first compute $x^{*}$ and and only then start the computation of the roots of $\phi$,
using $x^{*}$ as a starting point (assuming, as we did, that $\phi(x^{*}) < 0$).
The complexity of this detection phase is $\mathcal{O}(n)$.

We remark here that there exist other methods to detect if a polynomial has real roots, such as, e.g., Sturm sequences.
These methods have similar complexity, but they generally do not produce an appropriate starting point in $[r,R]$,
which is essential. 

\section{Computation of the roots}
\label{computation} 


In this section we assume that $\phi$ has two positive roots $r$ and $R$ and that we are given a point $\xb$
such that $r \leq \xb \leq R$, which can be obtained in the way explained at the end of the previous section.
We then propose a method to compute the roots
iteratively, with the iterates converging monotonically to the roots from inside the interval $[r,R]$.
The main idea behind this method is to approximate $\phi$ at a given point by a similar but simpler function that 
dominates $\phi$ on $[r,R]$. Its roots will be approximations to the roots of $\phi$, and the method then continues
iteratively from those approximations. The following theorem derives the approximation to $\phi$. 

\begin{theorem}
\label{phiapprox}
Let the polynomial     
\bdis
\phi(x)=x^{n} + \thnmo x^{n-1} + \dots + \thkpo x^{k+1} - \thk x^{k} + \thkmo x^{k-1} + \dots  + \thzero  \; , 
\edis
with $\thj \geq 0$ ($j=1,...,n$), $\thzero \thk \neq 0$, $n \geq 3$, and $1 \leq k \leq n-1$, have two
positive roots $r$ and $R$ with $r < R$ and let
\begin{eqnarray}
& & \phionex = x^{n} + \thnmo x^{n-1} + \dots + \thkpo x^{k+1} \; , \label{ph1def} \\ 
& & \phitwox = - \thk x^{k} + \thkmo x^{k-1} + \dots + \thone x + \thzero \; .  \label{phi2def} 
\end{eqnarray}   
Then for $r \leq \xb \leq R$, the trinomial   
\beq
\label{fdef}
f(x) = \alpha x^{n} - \beta x^{k} + \gamma \; , 
\eeq
where
\begin{eqnarray*}
& & \alpha = \dfrac{1}{n} \xb^{1-n} \phione'(\xb) > 0  \; , \\
& & \beta  = - \dfrac{1}{k} \xb^{1-k}\phitwo'(\xb) > 0  \; , \\
& & \gamma = \phi(\xb) - \xb \lb \frac{1}{n} \phione'(\xb) + \dfrac{1}{k} \phitwo'(\xb) \rb > 0 \; , 
\end{eqnarray*}
has two positive zeros $r_{1}$ and $r_{2}$ with $r \leq r_{1} < r_{2} \leq R$, and $f(x) \geq \phi(x)$
for $x \geq 0$.
\end{theorem}
\prf
Our goal is to approximate $\phi$ at a certain point $r \leq \xb \leq R$, for which $\phi(\xb) < 0$,
by a function $f$ that agrees at this point with $\phi$ in (at least) function and first derivative values. 
In addition, $f$ needs to dominate $\phi$ for $x \geq 0$. 
Clearly, a straightforward linear approximation is not possible as $\phi$ is composed of terms
that are both convex and concave. We solve this problem by constructing separate approximations for $\phione$
and $\phitwo$, whose sum is $\phi$, and which were defined in the statement of the theorem.   

The transformation of variables $w=x^{n}$ transforms $\phionex$ into 
\bdis
\phione(w^{1/n}) = w + \thnmo w^{\frac{n-1}{n}} + \dots + \thkpo w^{\frac{k+1}{n}} \; . 
\edis
This is a concave function of $w$, so that it is dominated by its linear approximation 
$\alpha_{1}w+\alpha_{2}$ at a point $\wb=\xb^{n}$.
A straightforward calculation shows that
\begin{eqnarray*}
& & \alpha_{1} = \lb \dfrac{1}{n} \wb^{\frac{1-n}{n}} \rb \phione'(\wb^{1/n}) = \dfrac{1}{n} \xb^{1-n} \phione'(\xb) \; , \\
& & \alpha_{2} = \phione(\wb^{1/n}) - \alpha_{1} \wb = \phionexb - \dfrac{1}{n} \xb \phione'(\xb)  \; ,
\end{eqnarray*}
where the derivatives are with respect to $x$.
Clearly, $\alpha_{1} > 0$ and
\bdis
\alpha_{2} = \phione(\xb) - \dfrac{1}{n} \xb \phione'(\xb)  
= \dfrac{1}{n} \lb \thnmo \xb^{n-1} + 2 \thnmt \xb^{n-2} + \dots + (n-k-1) \thkpo \xb^{k+1} \rb > 0 \; .
\edis
We have obtained an approximation to $\phionex$ of the form $\alpha_{1}x^{n}+\alpha_{2}$,
with $\alpha_{1}x^{n}+\alpha_{2} \geq \phi_{1}(x)$ for $x \geq 0$.

On the other hand, the transformation $y=x^{k}$ transforms $\phitwox$ into
\bdis
\phitwo(y^{1/k}) = - \thk y + \thkmo y^{\frac{k-1}{k}} + \dots + \thone y^{1/k} + \thzero \; ,
\edis
a concave function of $y$. It is therefore dominated by its linear approximation $\beta_{1}y+\beta_{2}$
at a point $\yb=\xb^{k}$, with 
\begin{eqnarray*}
& & \beta_{1} = \lb \dfrac{1}{k} \yb^{\frac{1-k}{k}} \rb \phitwo'(\yb^{1/k}) = \dfrac{1}{k} \xb^{1-k}\phitwo'(\xb) \; , \\
& & \beta_{2} = \phitwo(\yb^{1/k}) - \beta_{1} \yb = \phitwoxb - \dfrac{1}{k} \xb \phitwo'(\xb) \; . 
\end{eqnarray*}
Since
\bdis
\phitwo(\xb)- \dfrac{1}{k} \xb \phitwo'(\xb)  
= \dfrac{1}{k} \lb \thkmo \xb^{k-1} + 2\thkmt \xb^{k-2} + \dots + (k-1) \thone \xb + \thzero \rb \; ,
\edis
we have that $\beta_{2} = \phitwoxb - \dfrac{1}{k} \xb \phitwo'(\xb) > 0$,
from which $\dfrac{k}{\xb} \phitwoxb > \phitwo'(\xb)$. Because $\phixb = \phionexb + \phitwoxb \leq 0$, 
and therefore $\phitwoxb \leq -\phionexb < 0$, this means that $\phitwo'(\xb) < 0$, so that $\beta_{1} < 0$,
and we have obtained an approximation to $\phitwox$ of the form $\beta_{1}x^{k}+\beta_{2}$,
with $\beta_{1}x^{k}+\beta_{2} \geq \phi_{2}(x)$ for $x \geq 0$.
Consequently, our first-order approximation $f$ to $\phi$ at $\xb$ is given by the trinomial     
\bdis
f(x) = \alpha_{1} x^{n}+ \beta_{1} x^{k} + \alpha_{2} + \beta_{2}  \; ,
\edis
which corresponds to~(\ref{fdef}) in the statement of the theorem.
It satisfies $f(x) \geq \phix$, is of the same form as $\phi$, and, since $f(\xb) \leq 0$,
it has two roots $r_{1}$ and $r_{2}$ with $\rone \leq r_{1} < r_{2} \leq \rtwo$. \qed

An algorithm to compute $r$ and $R$ can now be based on Theorem~\ref{phiapprox} in a standard way:
starting from any initial point on the interval $(r,R)$, $f$'s largest root can be used as the next iterate 
in the computation of $R$. Convergence is monotonic from the left. An analogous algorithm is obtained for $r$
by considering the smallest root of $f$. In this case, convergence is monotonic from the right.
The function $\phi'$ has either one root, which must lie in $(r,R)$, or two roots, one in $(0,r)$ and one
in $(r,R)$, so that $\phi'(r) \phi'(R) \neq 0$. It is then a technical exercise to show that the order of
convergence of these algorithms is quadratic. 

However, to make the aforementioned algorithms implementable, we need a method to compute the roots of $f$ 
itself from inside the interval determined by its roots. The basis for such a method is provided 
by the following theorem, which derives a first order approximation to $f$ that dominates it, and whose roots
can be computed explicitly.

\begin{theorem}
\label{fapprox}
Let the trinomial $f(x) = \alpha x^{n} -\beta x^{k} + \gamma$ 
with $\alpha,\beta,\gamma > 0$, $n \geq 3$, and $1 \leq k \leq n-1$, have two
positive roots $r_{1}$ and $r_{2}$ with $r_{1} < r_{2}$.
Then for $r_{1} \leq \xb \leq r_{2}$, the function    
\beq
\label{hdef}
h(x) = \dfrac{\alpha \delta}{\epsilon - x^{k}} - \beta x^{k} + \gamma \; , 
\eeq
where
\begin{eqnarray*}
& & \delta = \dfrac{k}{n} \, \xb^{k+n} > 0  \; , \\
& & \epsilon = \dfrac{n+k}{n} \, \xb^{k} > \xb^{k}  \; , \\
\end{eqnarray*}
has two positive zeros $s_{1}$ and $s_{2}$ with $r_{1} \leq s_{1} < s_{2} \leq r_{2}$, and $h(x) \geq f(x)$
for $0 \leq x < (1+n/k)\xb$.
\end{theorem}
\prf
With the transformation of variables, $y=x^{k}$, $f$ can be written as $f(y^{1/k})=\alpha y^{n/k} - \beta y + \gamma$.
If the function $y^{n/k}$ is approximated at $\yb=\xb^{k}$ to first order by $\delta/(\epsilon-y)$, then this 
is equivalent to approximating $y^{-n/k}$ to first order by $(-1/\delta)y+\epsilon/\delta$. Since $y^{-n/k}$
is convex and its approximation is linear, this implies that   
$(-1/\delta)y+\epsilon/\delta \leq y^{-n/k}$, and therefore that $\delta/(\epsilon - y) \geq y^{n/k}$, or
\bdis
\dfrac{\delta}{\epsilon - x^{k}} \geq x^{n} \; .
\edis
Consequently, $h(x) \geq f(x)$ for $x \geq 0$ and, because $h(\xb)=f(\xb) \leq 0$, it has two positive 
roots $s_{1}$ and $s_{2}$
that satisfy $r_{1} \leq s_{1} < s_{2} \leq r_{2}$. The constants $\delta$ and $\epsilon$ are computed from 
the first order approximation conditions
\bdis
\dfrac{\delta}{\epsilon - \yb} = \yb^{n/k} 
\;\;\;\; \text{and} \;\;\;\;
\dfrac{\delta}{(\epsilon - \yb)^{2}} = \dfrac{n}{k} \yb^{n/k - 1} \; , 
\edis
as
\begin{eqnarray*}
& & \delta = \dfrac{k}{n} \, \yb^{1+n/k} = \dfrac{k}{n} \, \xb^{k+n} \; , \\ 
& & \epsilon = \dfrac{n+k}{n} \, \yb = \dfrac{n+k}{n} \, \xb^{k} \; .
\end{eqnarray*}
We note that $h(x)$ becomes unbounded as $x \longrightarrow (1+n/k)\xb$. This concludes the proof. \qed

\begin{figure}[H]
\begin{center}
\includegraphics[width=0.50\linewidth]{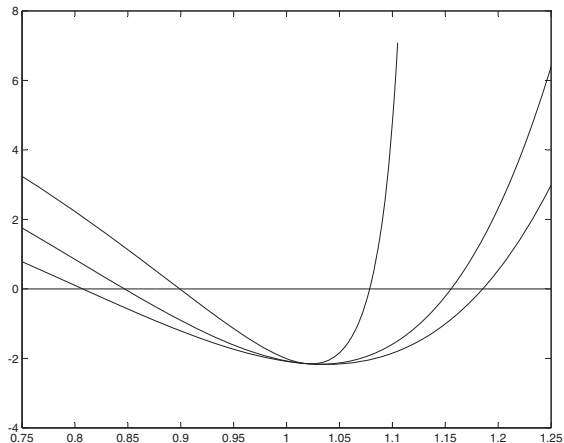}
\end{center}
\caption{The functions $\phi \leq f \leq h$ for $q(z)$, $k=3$, and $\xb=1.02$.}
\label{approximations}    
\end{figure}

Figure~\ref{approximations} shows the functions $\phi \leq f \leq h$ for the polynomial
\beq  
\label{qpolynomialdef}
q(z) = z^8 + z^7 + 3z^6 + \dfrac{1}{2} z^{4} +15z^{3} -2z^{2} +(i+1)z -4 \; ,
\eeq  
with $k=3$ and $\xb=1.02$.

The approximation $h$ of $f$ in Theorem~\ref{fapprox} leads to an iterative method for
the computation of the roots of $f$, by computing the roots of $h$ and then using the smallest and largest of 
those roots as the next iterates for the computation of $r_{1}$ and $r_{2}$, respectively.
As was the case for the roots of $\phi$, we obtain quadratic and monotonic convergence to the roots of $f$, 
although now the roots of the approximation can be computed explicitly since they are the roots of a quadratic.
This follows by setting $h(x)=0$, which is equivalent to
\bdis
\beta x^{2k} -(\epsilon \beta + \gamma) x^{k} + (\alpha \delta + \epsilon \gamma) = 0 \; ,
\edis
a quadratic equation in $x^{k}$. That both roots of this quadratic are real follows directly
from the properties of $h$.
Once available, the appropriate root of $f$ serves as the next iterate in the computation of the roots
of $\phi$.

The main computational effort in our approach to compute the positive roots of $\phi$ in Theorem~\ref{Pellet}
for a \emph{fixed} parameter $k$ is concentrated in the computation of the coefficients of the 
trinomial $f$, which requires $\mathcal{O}(n)$ arithmetic operations. The computation of the roots of $f$ 
itself is far less costly. 
Furthermore, the iterative process can be stopped at any moment, since each iterate provides a correct bound 
on the corresponding root of $\phi$. This is precisely what we set out to obtain.

$ $ \newline {\bf Example.}               

The motivation for this work was the absence of a method converging from inside the interval $[r,R]$,
which means that there is no equivalent method to compare our methods to. Instead we will illustrate
the method outlined above at the hand of an example, namely, the polynomial $q$ defined in~(\ref{qpolynomialdef}). 
For this polynomial, we computed the roots of the corresponding real polynomial $\phi$ for $k=3$ to a relative 
accuracy of $10^{-12}$. The steps involved in this process with the corresponding number of iterations are listed below. 
All computations are carried out to the same aforementioned relative accuracy of $10^{-12}$. 
\begin{tabbing}         
{\bf Step 1:} \=
Compute a starting point with Corollary~\ref{chibounds} for the computation of $x^{*}$, the positive \\ 
     \> root of $\chi$. \\
{\bf Step 2:} \>
Compute $x^{*}$ ($4$ Newton steps). \\ 
{\bf Step 3:} \>
Starting from $x^{*}$, compute the roots of $\phi$ with the help of the trinomials $f$, (6 iterations   \\
                        \> with the method based on Theorem~\ref{phiapprox} for each root). \\ 
\end{tabbing}        
\vskip -0.5cm
The number of iterations, necessary to compute the largest root of $f$ with the method from Theorem~\ref{fapprox}, 
was 9, 7, 5, 3, 2, and 1, corresponding to the six times such a root needed to be computed. For the smallest
root, we obtained 6, 5, 4, 2, 1, and 1 iterations. We observed that this decrease in the number of iterations 
required for the roots of $f$ is typical, regardless of the degree of the polynomial.

The number of iterations in the computation of the roots of $\phi$ does not vary significantly with increasing degree,
although the number of Newton steps necessary for the computation of $x^{*}$ tends to increase with increasing degree 
and may require an accelerated Newton method. On the other hand, it is usually not necessary to accurately 
compute $x^{*}$.

%
%




\end{document}